\documentclass[a4paper,10pt]{article}

\usepackage[notcite,notref]{showkeys}

\newtheorem{thm}{Theorem}

\newtheorem{rem}[thm]{Remark}
\newtheorem{lem}[thm]{Lemma}
\newtheorem{prop}[thm]{Proposition}
\def \a{{\alpha}}

\def \d{{\delta}}

\def \m{{\mu}}

\def \N{{\bf N}}

\def \qq{{\qquad}}

\def \noi{{\noindent}}
 \def \lan{{\langle}}
\def \ran{{\rangle}}

\def \a{{\alpha}}

\def \d{{\delta}}

\def \m{{\mu}}

\def \N{{\bf N}}

\def \qq{{\qquad}}

\def \noi {{\noindent}}

\def\qed{\hbox{\vrule height 6pt depth 0pt width 6pt}}
\def\cqfd{\hfill\penalty 500\kern 10pt\qed\medbreak}

\font\phh=cmcsc10  at  8 pt

  at 10,5 pt

 \font\phh=cmcsc10
at  8 pt

\title{On an identity of Ky Fan}
\author{Michel Weber}

\begin{document}

\maketitle
%%%%%%%%%%%%%%%%%%%%%%%%%%%%%%%%%%%%%%%%%%%%%%%%%%%%%%%%%%%%%%%%%%%%%%%%%%%%%%%%%%%%%%%%%%%%%%

\begin{abstract}
We give several applications of an identity  
 for sums of weakly stationary sequences due to Ky Fan.
\end{abstract}

%%%%%%%%%%%%%%%%%%%%%%%%%%%%%%%%%%%%%%%%%%%%%%%%%%%%%%%%%%%%%%%%%%%%%%%%%%%%%%%%%%%%%%%%%%%%%%%

\section{Introduction and  results}
\scrollmode
 \medskip \noindent
  Let $X_1, X_2, \ldots$ be a weakly stationary sequence in an Hilbert $H$. In   \cite{KF2} (see p.598), Ky Fan
observed that 
    for any  two   positive integers   
$n,m$, one has
$$ {\big\| X_1+\ldots + X_{n }\big\|^2\over n}+{\big\| X_1+\ldots + X_{m }\big\|^2\over m} -{\big\| X_1+\ldots + X_{n+m
}\big\|^2\over n+m}\qq\qq\qq$$
\begin{equation}\qq\qq\qq= {n(n+m)\over m}\ \Big\|{  X_1+\ldots + X_{n } \over n} -{  X_1+\ldots + X_{n+m } \over
n+m}\Big\|^2 . 
\end{equation}
 
\noi This nice identity was   applied in the same paper. Rather surprisingly, this result did not seem   to have caught much attention.
 
The object of this short note is to indicate some other consequences of identity (1), which have not been quoted in \cite{KF1} or
\cite{KF2}. Before going further, and since no proof of this identity is given in \cite{KF2}, we thought worth to give one. The proof
goes as follows. Put for any positive integer
$n$,
$S_n =X_1+\ldots + X_n$, and if $m$ is another positive integer let $T_{n,m}= S_{n+m}-S_n$, so that $S_{n+m}= S_n+ T_{n,m}$. Then
$$\big\| {S_n\over n}-{S_{n+m}\over n+m}\big\|^2=  {\|S_n\|^2\over n^2} + {\|S_{n+m}\|^2\over (n+m)^2}- \lan {S_n\over n},
{S_{n+m}\over n+m}\ran -\lan {S_{n+m}\over n+m} , {S_n\over n}\ran ,$$
and so
$${n(n+m)\over m}\ \big\| {S_n\over n}-{S_{n+m}\over n+m}\big\|^2\qq\qq\qq\qq\qq\qq\qq\qq\qq\qq\qq $$
\vskip -17pt\begin{eqnarray*} & =&
  {(n+m) \|S_n\|^2\over n m}  + { n\|S_{n+m}\|^2\over
m(n+m) }    -{1\over m}\big\{\lan  S_n ,  S_{n+m} \ran +\lan  S_{n+m}  ,  S_n \ran \big\}\cr 
& =& {  \|S_n\|^2\over n  }+ \Big({\|S_n\|^2\over m  }- {\|S_m\|^2\over m  }\Big)+ {\|S_m\|^2\over m  } -{  \|S_{n+m}\|^2\over  n+m 
} 
\cr 
& \quad  +& {  \|S_{n+m}\|^2\over   m } -{1\over m}\Big\{\lan  S_n ,  S_{n+m} \ran +\lan  S_{n+m}  ,  S_n \ran \Big\}.
\end{eqnarray*}
But 
$$\lan  S_n ,  S_{n+m} \ran +\lan  S_{n+m}  ,  S_n \ran =  2\|S_n\|^2 + \lan  S_n ,  T_{n,m} \ran    + \lan     T_{n,m}, 
S_n 
\ran , $$
so that in turn
\begin{eqnarray}{n(n+m)\over m}\ \big\| {S_n\over n}-{S_{n+m}\over n+m}\big\|^2  
& =& {  \|S_n\|^2\over n  }+   {\|S_m\|^2\over m  } -{  \|S_{n+m}\|^2\over 
n+m  }\cr 
& \ \ &  \quad + {1\over m}\Big(\|S_{n+m}\|^2 -\|S_{n }\|^2 
\cr 
& \ \ &  \quad-  \|S_{ m}\|^2-  \lan  S_n ,  T_{n,m} \ran    - \lan     T_{n,m}, S_n 
\ran\Big) \cr 
& =& {  \|S_n\|^2\over n  }+   {\|S_m\|^2\over m  } -{  \|S_{n+m}\|^2\over 
n+m  } ,  
\end{eqnarray}
since  $S_{n+m}= S_n+ T_{n,m}$. And we are done.

\medskip Note that the weak stationarity assumption was only used in the last line of
calculations, to say that $\|T_{n,m}\|=\|S_m\|$. And consequently, if $X_1, X_2, \ldots$ is any sequence in  
$H$ satisfying      
 \begin{equation} \|\sum_{k= 1}^{ m}X_k\|\ge \|\sum_{k=n+1}^{n+m}X_k\|,\qq (n\ge 1, m\ge 1),
\end{equation}
then, for any      positive integers $n,m$

\begin{equation}   {n(n+m)\over m}\ \Big\|{  S_n \over n} -{  S_{n+m } \over
n+m}\Big\|^2 \le { \| S_n\|^2\over n}+{ \| S_m\|^2\over m} -{ \| S_{n+m }\|^2\over
n+m}  .
\end{equation}
\medskip\par
 So is the case when for instance $X_i=T^i X_0$, $i=1, \ldots$ where $T$ is a contraction in $H$, and $X_0$ some fixed element of $H$.  A 
natural question concerns the possibility to replace the  norming  factor
$n^{-1 }$ by another one. The lemma below shows that (4) remains true with  norming factor $n^{-\d}$,   $1/2<\d\le 1$.
\begin{lem}   Let $\{\a_k, k\ge 1\}$ be a sequence of positive reals satisfying the following condition:
\begin{equation} \a_{n+m}\le \a_{n }+\a_{ m},\qq n,m\ge 1. 
\end{equation}
  Let $X_1, X_2, \ldots$ be a sequence in  
$H$ satisfying  (3). Then, for any      positive integers $n,m$ 
\begin{equation}   {\a_{n }\a_{n+m}\over \a_{ m}}\, \Big\|{  S_n \over \a_{n }} -{  S_{n+m } \over
\a_{n+m}}\Big\|^2 \le { \| S_n\|^2\over \a_{n }}+{ \| S_m\|^2\over \a_{ m}} - { \| S_{n+m }\|^2\over
\a_{n+m}} \left({\a_{n+m}-\a_{n }\over \a_{ m}}\right)  . 
\end{equation} 
\end{lem}
\medskip \begin{rem}--- \rm A typical case where Lemma 1 applies is when $\a_n=n^\d$ with $0\le \d\le 1$, and we get in
particular for any
$m\ge n\ge 1$
\begin{equation}   {n^\d (n+m)^\d\over m^\d}\, \Big\|{  S_n \over n^\d} -{  S_{n+m } \over
(n+m)^\d}\Big\|^2 \le { \| S_n\|^2\over n^\d}+{ \| S_m\|^2\over m^\d} -  { \| S_{n+m }\|^2\over
(n+m)^\d} \left({ (n+m)^\d - n^\d\over   m ^\d}\right)   . 
\end{equation}  
--- Notice also that (7) with $\d=1$ is (4).
\end{rem}
 \medskip
\noi {\it Proof.} Similarly as before 
  \begin{eqnarray} {\a_{n }\a_{n+m}\over \a_{ m}}\, \big\| {S_n\over \a_{ n}}-{S_{n+m}\over \a_{n+m}}\big\|^2&   
  =& {  \|S_n\|^2\over \a_{n }  } + {\|S_m\|^2\over\a_{ m}  } -{ 
\|S_{n+m}\|^2\over  \a_{n+m} }\cr & &\qq +{ 
\|S_{n+m}\|^2\over   \a_{ m}\a_{n+m} } \Big(\a_n-( \a_{n+m}-\a_{ m} )\Big)
\cr 
&  &  +   \|S_{n }\|^2\Big({\a_{n+m}-(\a_{n }+\a_{ m})\over \a_{n }\a_{ m} }\Big)\cr & &\qq +{1\over \a_{ m} }\big(\|T_{n,m}\|^2-
\|S_{ m}\|^2
\big).
\end{eqnarray}
Using assumptions (3) and (5) gives
 $${\a_{n }\a_{n+m}\over \a_{ m}}\, \big\| {S_n\over \a_{ n}}-{S_{n+m}\over \a_{n+m}}\big\|^2 \le {  \|S_n\|^2\over \a_{n }  } + 
{\|S_m\|^2\over\a_{ m}  }  - { \| S_{n+m }\|^2\over
\a_{n+m}} \left({\a_{n+m}-\a_{n }\over \a_{ m}}\right),$$  as claimed.\cqfd 
\smallskip
 A simple although quite interesting consequence of Ky Fan's identity is  
\begin{equation}{  \|S_{n+m}\|^2\over 
n+m  }\le {  \|S_{n }\|^2\over 
n   }+ {  \|S_{ m}\|^2\over 
 m  } ,\end{equation}
which is valid for any two positive integers $n,m$. This is inequality (4.8) in \cite{KF2}. 
\smallskip
Recall that a sequence    $\{g_n ,n\ge 1\}$ of
real numbers is  subadditive when
$$ g_{n+m}\le g_{n }+ g_{ m}.  $$
Then  we have the well-known lemma
\begin{lem}If $\{g_n ,n\ge 1\}$ is a subadditive sequence of real numbers, then $g_n/n$ converges to
$\displaystyle{\inf_{n\ge 1}(g_n/n)}$.  
\end{lem}
\noi {\it Proof.} Fix an arbitrary positive integer $N$ and write $n= j_n N + r_n$ with $1\le r_n\le N$. Clearly ${j_n\over n}\to 
{1\over N}$ as
$n$ tends to infinity. Further
$$\inf_{n\ge 1}{g_n\over n }\le {g_n\over n }\le {g_{j_n N}+ g_{r_n  }\over n } \le {g_{j_n N}\over j_n N } + {g_{r_n  }\over  n   }
\le {j_ng_{  N}\over j_n N } + {g_{r_n  }\over  n   } ={ g_{  N}\over   N } + {g_{r_n  }\over  n   }. $$
Letting now $n$ tend to infinity gives 
 $\displaystyle{\inf_{n\ge 1}{g_n\over n }\le \liminf_{n\to \infty}{g_n\over n }\le \limsup_{n\to \infty}{g_n\over n } \le { g_{ 
N}\over   N }  .  }$ As $N$ was arbitrary, the lemma is proved.\cqfd
\medskip We thus deduce from   (2) and this lemma  applied to
$g_n:={ 
\|S_{n }\|^2\over  n   }$  that 
 \begin{equation} \lim_{n\to \infty} \|{S_n\over n}\|= \inf_{n\ge 1} \|{S_n\over n}\|.  
 \end{equation}
   This is a remarkable direct consequence of Ky Fan's identity, which  remains true for averages of contractions.  If $T$ is a 
contraction in $H$, in  view
of   the Riesz's decomposition (\cite{Kr}, lemma 1.3 p.4), the orthogonal complement $H_T^\perp$ of  $  H_T =\{g\in H: Tg=g\}$ coincides
with the closure of the subspace   spanned by  
$ \{h-Th: h\in H\}$. Then it suffices to proceed by approximation.
\medskip
 In the next proposition we examine the ratios
$${ \big\|
{S_{n_{k}}\over n_{k}}-{S_{n_{k+1}}\over n_{k+1}}\big\|^2  \over \big({1\over n_{k }}- {1\over n_{k+1} }\big)},$$
 where ${\cal N}=\{n_{k}, k\ge 1\}$ is an  increasing sequence of positive integers. Notice that in the orthonormal case, namely if  $X_1,
X_2, \ldots$  is an orthonormal sequence, then precisely
$\big\| {S_{n_{k}}\over n_{k}}-{S_{n_{k+1}}\over n_{k+1}}\big\|^2={1\over n_{k }}- {1\over
n_{k+1} } $.

\begin{prop}  (a) Let
${\cal N}=\{n_{k}, k\ge 1\}$ be an  increasing sequence of positive integers. Then
$$ \limsup_{N\to \infty} {  1\over n_{N}  } \sum_{k=1}^{N -1} { \big\|
{S_{n_{k}}\over n_{k}}-{S_{n_{k+1}}\over n_{k+1}}\big\|^2  \over \big({1\over n_{k }}- {1\over n_{k+1} }\big)}
   \le    \limsup_{N\to \infty}   \sup_{ 1\le k<N } \Big| {\|S_{n_{k+1}- n_{k}}\|^2 \over (n_{k+1}- n_{k})^2  } -  { 
\|S_{n_{N}}\|^2\over n^2_{N}  }\Big|  .
  $$
(b) Further, if $\lim_{k\to \infty} n_{k+1}-n_k=\infty$,   
$$\lim_{N\to \infty}
 {  1\over n_{N}  } \sum_{k=1}^{N -1} { \big\|
{S_{n_{k}}\over n_{k}}-{S_{n_{k+1}}\over n_{k+1}}\big\|^2  \over \big({1\over n_{k }}- {1\over n_{k+1} }\big)}
    =0.  $$
(c) And  
$$  \lim_{N,a\to \infty} 
   {  1\over Na  } \sum_{k=1}^{N  -1}  \big\|
 k( S_{a(k+1)}-S_{ak})-S_{ak}\big\|^2   
    =0  .
$$
(d) Finally  let ${\cal D}=\big\{D_j, j\ge 1\big\}$ be a chain: $D_j|D_{j+1}$ for every $j$. Then 
$$ \sum_{j=1}^\infty{  1\over D_{j+1} 
}
\sum_{k=1}^{{D_{j+1}\over D_j} -1} { \big\| {S_{kD_j}\over kD_j}-{S_{(k+1)D_j}\over(k+1)D_j}\big\|^2  \over {1\over kD_j}- {1\over (k+1)D_j }}= \big({\|S_{
D_1}\|  \over D_1   }  \big)^2-\lim_{J\to
\infty}
\big({ 
\|S_{D_{J+1}}\| \over 
 D_{J+1}    } \big)^2<\infty. $$ 
\end{prop} 
\begin{rem} \rm --- A sequence $\{a_n,n\ge 1\}$ of real numbers   converges  in      
density  to 0, which we write $ {D-\lim_{n\rightarrow \infty} a_n=0}$, if there exists a subset ${\cal J}$ of $\N$ of density
one, such that 
 $ { \lim_{{\cal J}\ni
n\to \infty} a_n = 0} $. For {\it bounded} sequences, it is an exercice to show that
   $\displaystyle{D-\lim_{n\rightarrow \infty} a_n=0}$,
if  and only if,
 $ \lim_{n\rightarrow
\infty} {1\over n}\sum_{j=0}^{n-1}|a_n |=0    $. Since 
$$      { \big\|
  {S_{ak}\over ak}-{S_{a(k+1)}\over a(k+1)}\big\|^2  \over  {1\over ak}- {1\over a(k+1) } }=  \big\|
 k( S_{a(k+1)}-S_{ak})-S_{ak}\big\|^2 ,$$    
  $(c)$    means
$$  \lim_{N,a\to \infty} 
   {  1\over Na  } \sum_{k=1}^{N  -1}     { \big\|
  {S_{ak}\over ak}-{S_{a(k+1)}\over a(k+1)}\big\|^2  \over  {1\over ak}- {1\over a(k+1) } }   
    =0  .
$$
Thus  along   linearly
growing sequences, the averages of weakly stationary sequences   have, in density, increments comparable to averages of
orthogonal sequences, which is a bit unexpected.  
\end{rem} 
\noi {\it Proof.} ---   From  Ky Fan's identity, we   get for each
$k$  
 $$   { \big\| {S_{n_{k}}\over n_{k}}-{S_{n_{k+1}}\over n_{k+1}}\big\|^2  \over {1\over n_{k }}- {1\over n_{k+1} }}
   = \left({  \|S_{n_{k}}\|^2\over n_{k}  }- {  \|S_{n_{k+1}}\|^2\over 
n_{k+1}  }\right) +   {\|S_{n_{k+1}- n_{k}}\|^2\over n_{k+1}- n_{k}  }   . $$
Summing   from $k=1$ up to $N-1$ leads to 
\begin{equation}  \sum_{k=1}^{N -1} { \big\| {S_{n_{k}}\over n_{k}}-{S_{n_{k+1}}\over n_{k+1}}\big\|^2  \over {1\over n_{k }}- {1\over
n_{k+1} }}
    =     {  \|S_{n_{1}}\|^2\over n_{1}  }- {  \|S_{n_{N}}\|^2\over 
n_{N}  }  +  \sum_{k=1}^{N -1} (n_{k+1}- n_{k}) \Big({\|S_{n_{k+1}- n_{k}}\| \over n_{k+1}- n_{k}  }\Big)^2  
\end{equation} 
Dividing both sides by $n_{N}$      gives 
\begin{eqnarray} {  1\over n_{N}  } \sum_{k=1}^{N -1} { \big\|
{S_{n_{k}}\over n_{k}}-{S_{n_{k+1}}\over n_{k+1}}\big\|^2  \over {1\over n_{k }}- {1\over n_{k+1} }}
 &  =    & {  \|S_{n_{1}}\|^2\over n_{1}  n_{N} }- {  \|S_{n_{N}}\|^2\over 
n^2_{N}  }  \cr & & \qq+ {  1\over n_{N}  } \sum_{k=1}^{N -1} (n_{k+1}- n_{k}) \Big({\|S_{n_{k+1}- n_{k}}\| \over n_{k+1}- n_{k} 
}\Big)^2 
\cr 
  & \le & {  \|S_{n_{1}}\|^2\over n_{1} n_{N}  }      + {  1\over n_{N}  } \sum_{k=1}^{N -1} (n_{k+1}-
n_{k}) \bigg\{\Big({\|S_{n_{k+1}- n_{k}}\| \over n_{k+1}- n_{k}  }\Big)^2    \cr & & \qq -{  \|S_{n_{N}}\|^2\over n^2_{N}  }\bigg\}  .
 \end{eqnarray} 
Letting next $N$ tend to infinity yields 
$$\limsup_{N\to \infty} {  1\over n_{N}  } \sum_{k=1}^{N -1} { \big\|
{S_{n_{k}}\over n_{k}}-{S_{n_{k+1}}\over n_{k+1}}\big\|^2  \over \big({1\over n_{k }}- {1\over n_{k+1}
}\big)}\qq\qq\qq\qq\qq\qq\qq\qq\qq\qq
$$
\vskip -17pt \begin{eqnarray} 
 &  \le  & \limsup_{N\to \infty}  {  1\over n_{N}  } \sum_{k=1}^{N -1} (n_{k+1}-
n_{k}) \bigg\{\Big({\|S_{n_{k+1}- n_{k}}\| \over n_{k+1}- n_{k}  }\Big)^2-  {  \|S_{n_{N}}\|^2\over n^2_{N}  }\bigg\}  .
  \cr 
&  \le  &  \limsup_{N\to \infty}   \sup_{ 1\le k<N } \Big| {\|S_{n_{k+1}- n_{k}}\|^2 \over (n_{k+1}- n_{k})^2  } -  { 
\|S_{n_{N}}\|^2\over n^2_{N}  }\Big|  .
 \end{eqnarray} 
Hence the first claim is proved. 

\smallskip\noi  --- If $\lim_{k\to \infty} n_{k+1}-n_k=\infty$ and suppose first that $\lim_{n\to \infty} {\|S_n \| \over  
n   }=0$. Then  
 $\lim_{k\to \infty} {\|S_{n_{k+1}- n_{k}}\| \over n_{k+1}- n_{k}  } =
0, $ and so letting $N$ tend to infinity in the first equality in (11) gives
$$\lim_{N\to \infty} {  1\over n_{N}  } \sum_{k=1}^{N -1} { \big\|
{S_{n_{k}}\over n_{k}}-{S_{n_{k+1}}\over n_{k+1}}\big\|^2  \over \big({1\over n_{k }}- {1\over n_{k+1} }\big)}
=0  .
$$
Hence the second claim of the proposition is proved in that case. If $\lim_{n\to \infty} {\|S_n \| \over   n   }>0$, there exists
$\chi\in H$ such that $\lim_{n\to \infty} \|{ S_n \over   n   }-\chi \| =0$. Indeed, first recall (\cite{Kr}, p.32) that $\{X_i , i\ge 1\}$ may be represented
as a sequence
$\{T^iX_1, i\ge 0\}$ where $T$ is an isometry in some Hilbert space, which we denote again $H$.  Next  by  the mean ergodic
theorem of von Neumann (\cite{Kr}, p.4),  the limit
$\chi$ is identified as the projection on the subspace $H_T=\{g\in H: Tg=g\}$ of $X_1$. Applying the result previously obtained to the
weakly stationary sequence
$\{X_i-\chi, i\ge 1\}$, allows to reach the same conclusion in this case as well.  

\smallskip\noi  --- Now assume   $n_k= ak$, $a$ being some fixed positive integer. Replace $n_k$ by its value in the first part of (11). 
\begin{eqnarray} {  1\over Na  } \sum_{k=1}^{N -1} { \big\|
{S_{ ka}\over  ka}-{S_{ (k+1)a}\over  (k+1)a}\big\|^2  \over \big({1\over  ka}- {1\over  (k+1)a }\big)}
 &  =   &  {  \|S_{a}\|^2\over a^2{N} }- {  \|S_{a{N}}\|^2\over 
a^2 {N}^2  }  + {  1\over   {N}  } \sum_{k=1}^{N -1}    {\|S_{a}\|^2 \over a^2  } \cr  &  =   &  {\|S_{a}\|^2 \over a^2  } -{ 
\|S_{a{N}}\|^2\over  a^2 {N}^2  } 
\end{eqnarray} Hence 
 
$$ \limsup_{N,a\to \infty } {  1\over Na  } \sum_{k=1}^{N -1} { \big\|
{S_{ ka}\over  ka}-{S_{ (k+1)a}\over  (k+1)a}\big\|^2  \over \big({1\over  ka}- {1\over  (k+1)a }\big)}
    \le    \limsup_{N,a\to \infty }     \Big| {\|S_{a}\|^2 \over a^2  } -  { 
\|S_{ Na}\|^2\over  ( Na)^2  }\Big|    =0  .
$$
The expression in the right-hand side  being also rewritten as
$$ {  1\over Na  } \sum_{k=1}^{N  -1}  \big\|
 k( S_{a(k+1)}-S_{ak})-S_{ak}\big\|^2  ,$$
we get (c).
\smallskip\par \noi ---
Let $\{D_j, j\ge 1\}$ be a chain, and applies   (12) with $N={D_{j+1}/D_j}$, $a=D_j$.
\begin{equation} {  1\over Na  } \sum_{k=1}^{N -1} { \big\|
{S_{ ka}\over  ka}-{S_{ (k+1)a}\over  (k+1)a}\big\|^2  \over \big({1\over  ka}- {1\over  (k+1)a }\big)}
      =  {\|S_{a}\|^2 \over a^2  } -{  \|S_{a{N}}\|^2\over 
a^2 {N}^2  } 
\end{equation}  We obtain
$$  {  1\over D_{j+1}  } \sum_{k=1}^{{D_{j+1}\over D_j} -1} { \big\|
{S_{kD_j}\over kD_j}-{S_{(k+1)D_j}\over(k+1)D_j}\big\|^2  \over {1\over kD_j}- {1\over (k+1)D_j }}
    =     \big({\|S_{ D_j}\|  \over D_j   }  \big)^2  - \big({  \|S_{D_{j+1}}\| \over 
 D_{j+1}    } \big)^2     . 
    $$ 
Summing up from $j=1$ to $j=J$ gives
$$ \sum_{j=1}^J{  1\over D_{j+1}  } \sum_{k=1}^{{D_{j+1}\over D_j} -1} { \big\|
{S_{kD_j}\over kD_j}-{S_{(k+1)D_j}\over(k+1)D_j}\big\|^2  \over {1\over kD_j}- {1\over (k+1)D_j }}
=   \big({\|S_{ D_1}\|  \over D_1   }  \big)^2- \big({  \|S_{D_{J+1}}\| \over 
 D_{J+1}    } \big)^2  .$$
 Therefore
$$ \sum_{j=1}^\infty{  1\over D_{j+1}  } \sum_{k=1}^{{D_{j+1}\over D_j} -1} { \big\|
{S_{kD_j}\over kD_j}-{S_{(k+1)D_j}\over(k+1)D_j}\big\|^2  \over {1\over kD_j}- {1\over (k+1)D_j }}= \big({\|S_{ D_1}\|  \over D_1   }  \big)^2-\lim_{J\to
\infty}
\big({ 
\|S_{D_{J+1}}\| \over 
 D_{J+1}    } \big)^2<\infty. $$ 
Hence (d), the proof is now complete.  \cqfd 

For arithmetic progressions, Proposition 4 shows that 
 $$   { \big\|
{S_{ ka}\over  ka}-{S_{ (k+1)a}\over  (k+1)a}\big\|^2  \over \big({1\over  ka}- {1\over  (k+1)a }\big)}   $$
"asymptotically" converges in density to $0$.
The question naturally arises when for {\it any} increasing
sequence ${\cal N}$, a convergence in density to $0$ do hold.
   \bigskip

In our next result, we give an example having this property. Recall   that a   Dunford-Schwartz contraction  is a linear contraction $T$ on
$L_1$ of a $\sigma$-finite measure space, with
$\|Tf\|_\infty \le \|f\|_\infty$ for $f \in L_1 \cap L_\infty$
, and induces
a contraction on all $L_p$, $1<p \le \infty$ (\cite{DS};
 see also \cite{Kr}  p. 65  for $T$ positive).  
The limit $E(T)f:= \lim_n {1\over n} \sum_{k=1}^n T^k f$  exists
a.e. for $f \in L_p$, $1 \le p < \infty$, and also in $L_p$-norm for
$p>1$ (and in $L_1$-norm in probability spaces). Write $S_n^T=\sum_{l=1}^n T^l$. Let $0<\a<1$. 
By  Corollary 2.15 in \cite{DL}, when $T$ is induced on $L_p$
by a Dunford-Schwartz operator, if $f \in (I-T)^\alpha L_p$,  then    
$$ \lim_{n\to \infty} \|{S_n^T f\over n^{1-\alpha}}  \|_p = 0 .$$
 
 If $T$ is an isometry on $L_2$ of a $\sigma$-finite measure space, for instance if $T$ is the isometry induced by an ergodic
transformation, 
then $\sqrt{I-T}L_2$ is a dense sub-space in the space of elements $f\in L_2$ such that $\lan f,1\ran=0$; and strictly contains 
the space of coboundaries. And we have the following characterization  (\cite{C} Proposition 2.2): $f\in
\sqrt{I-T}L_2$ if and only if the series
$\sum_{n=1}^\infty \|S_n^T f\|^2/n^2 $ converges.   A spectral
characterization is further given in   Proposition 2.2 of \cite{C} (see also  \cite{DL} Theorem 4.4). 
\begin{thm}  Let $T$ be an isometry on $L_2$ of a $\sigma$-finite measure space, and $f \in \sqrt{(I-T)} L_2$. Then for any increasing
sequence ${\cal N}=\{n_{k}, k\ge 1\}$   of positive integers, we have 
$$\lim_{K\to \infty}
{1\over   K  } \sum_{k=1}^{K }  {\displaystyle{ \Big\| {S^T_{n_{k}}f\over n_{k}}-{S^T_{n_{k+1}}f\over n_{k+1}}\Big\|^2 } \over
\displaystyle{
\Big({1\over n_{k }}- {1\over n_{k+1} }\Big)} }    =0.   $$
\end{thm} 
\noi {\it Proof.}  We apply the previous remark with $p=2$, $\a=1/2$. By (11) 
\begin{eqnarray*} \sum_{k=1}^{K-1} { \big\| {S_{n_{k}}\over n_{k}}-{S_{n_{k+1}}\over n_{k+1}}\big\|^2  \over {1\over n_{k }}- {1\over
n_{k+1} }}
 &  =   &  {  \|S_{n_{1}}\|^2\over n_{1}  }- {  \|S_{n_{K}}\|^2\over 
n_{K}  }  +  \sum_{k=1}^{K -1}   {\|S_{n_{k+1}- n_{k}}\|^2 \over n_{k+1}- n_{k}  } \cr &\le  &   {  \|S_{n_{1}}\|^2\over n_{1}  } 
+o(K)   
  .
\end{eqnarray*}  
Dividing both sides by $K$ and letting next $K$ tend to infinity  achieves the proof.
\cqfd

\bigskip
\bigskip    We conclude with another inequality. Put
$$ f^2(n)= \sup_{n'\le n} {\|S_{n'}\|^2\over n'}.$$
\begin{lem}   $f^2$ is subadditive: we have
 $ f^2(n+m)\le f^2(n)+f^2(m)  $  for any $m,n\ge 1$.
\end{lem}
\noi {\it Proof.}  There is no loss to assume $n\le m$. Let $\m\le m+n$. Consider three cases:
\smallskip
\smallskip
$i)$ $\m\le n$. Then $\|S_\m\|^2/\m\le f(n)\le f(n)+ f(m)$.
\smallskip

$ii)$ $n<\m\le m$. Write $\m=\m-n+n$. We have $\m-n\le m$ and using (9)
$${\|S_\m\|^2\over \m}={\|S_{(\m-n)+n}\|^2\over (\m-n)+n}\le {\|S_{ \m-n  }\|^2\over  \m-n  } + {\|S_n\|^2\over n}\le f^2(\m-n)+ f^2(n)\le
f^2(m)+f^2(n).$$

  $iii)$ $m<\m\le n+m$. Then $\m=\m-m+m$. We have $\m-m\le n+m-m=n$ and using (9) again
$${\|S_\m\|^2\over \m} \le {\|S_{ \m-m  }\|^2\over  \m-m  } + {\|S_m\|^2\over m}\le f^2( \m-m)+ f^2(m)\le
f^2(n)+f^2(m).$$

\cqfd
\medskip
  Put 
$${\rm  I}(f,x,y)= {f^2(x+y)-f^2(x )+f^2( y)\over 2f(x+y)f(y)}.$$
Observe that if $h (x)=\  |x|^{1/2}$, then ${\rm  I}(h,x,y)=2\sqrt {{y\over y+x}}$. 

\begin{prop} For any positive integers $x,y$ such that $f^2(x)/x\ge f^2(y)/y$, we have 
$${\rm  I}(f,x,y)\le 2\sqrt {{y\over y+x}}.$$
\end{prop}
\medskip
\noi {\it Proof.}   
Write 
$$f^2(x+y)=p(x+y), \qq f^2(y)= qy, \qq f^2(x)=rx. $$
Substituing these values into ${\rm  I}(f,x,y)$ gives  
$${\rm  I}(f,x,y)={ph^2(x+y)-rh^2(x )+qh^2( y)\over 2\sqrt{pq}h(x+y)h(y)} $$ 
Consider the following expression 
$$J=\Big({\rm  I}(f,x,y)- {\rm  I}(h,x,y)\Big) 2 h(x+y)h(y)=p(x+y)+qy -rx-2y\sqrt{pq}$$
It suffices to prove $J\le 0$. By lemma 3, $p(x+y) \le rx +qy$.  

\smallskip
\noi --- If $q\le p$, then $J\le 2qy -2y\sqrt{pq}= 2y\sqrt q\big[\sqrt q -\sqrt p\big]\le 0$.
\smallskip
\noi--- If $q\ge p$, then
$$ J =-\sqrt p(\sqrt q-\sqrt p)(x+y)+ \sqrt q(\sqrt q-\sqrt p)y+ (\sqrt{pq}-r)x.$$
We have $\sqrt q- \sqrt p\ge 0$, moreover $qy= f^2(y)\le f^2(x+y)=p(x+y)$. Thus 
\begin{eqnarray*} J &\le &\big[-\sqrt p(\sqrt q-\sqrt p){q\over p}+ \sqrt q(\sqrt q-\sqrt p)\big]y+ (\sqrt{pq}-r)x \cr &=&-{\sqrt q\over
\sqrt p}\big(\sqrt q -\sqrt p\big)^2y+ (\sqrt{pq}-r)x.
\end{eqnarray*} 
But we assumed $f^2(x)/x\ge f^2(y)/y$, thus $r\ge q\ge p$. And so $\sqrt{pq}-r\le 0$. Therefore $J\le 0$ as required.\cqfd
\bigskip\noi {\it Acknowledgments.} I thank Christophe Cuny and Michael Lin  for   useful comments. 

%%%%%%%%%%%%%%%%%%%%%%%%%%%%%%%%%%%%%%%%%%%%%%%%%%%%%%%%%%%%%%%%%%%%%%%%%%%%%%%%%%%%%%%%%%%%%%%%
%%%%%%%%%%%%%%%%%%%%%%%%%%%%%%%%%%%%%%%%%%%%%%%%%%%%%%%%%%%%%%%%%%%%%%%%%%%%%%%%%%%%%%%%%%%%%%%%
%%%%%%%%%%%%%%%%%%%%%%%%%%%%%%%%%%%%%%%%%%%%%%%%%%%%%%%%%%%%%%%%%%%%%%%%%%%%%%%%%%%%%%%%%%%%%%%%

{
}
\bigskip
  \noi {\phh Michel  Weber, \noi  Math\'ematique (IRMA),
Universit\'e Louis-Pasteur et C.N.R.S.,   7  rue Ren\'e Descartes, 
67084 Strasbourg Cedex, France. 
\par\noindent
E-mail: \  \tt weber@math.u-strasbg.fr}
\end{document}